\documentclass[12pt,a4paper]{article}
\usepackage[utf8]{inputenc}
\usepackage{amsmath}
\usepackage{amsfonts}
\usepackage{amssymb}
\usepackage{url}
\usepackage{graphicx}
\graphicspath{ {images/} }

\newtheorem{remark}{Remark}
\newtheorem{lemma}{Lemma}
\newtheorem{corollary}{Corollary}
\newtheorem{theorem}{Theorem}
\newtheorem{example}{Example}
\newtheorem{algorithm}{Algorithm}

\author{Alexander Kushkuley \\ (Salesforce/Demandware, akushkuley@salesforce.com)}
\title{Heavy Hitters and Bernoulli Convolutions  }
\begin{document}
\maketitle

\begin{abstract}
\noindent	A very simple event frequency approximation algorithm that is sensitive to event timeliness is suggested. The algorithm iteratively updates categorical click-distribution, producing (path of) a random walk on a standard $n$-dimensional simplex. Under certain conditions, this random walk is self-similar and corresponds to a biased Bernoulli convolution. Algorithm evaluation naturally leads to estimation of moments of biased (finite and infinite) Bernoulli convolutions.
\end{abstract}

\section{Introduction}
To quote \cite{Freq2}, "there is a need to estimate the count of a given item $i$ (or event or combination thereof) during some period of time $t$...Typically, items with highest counts, commonly known as heavy hitters, are of most interest".
\newline
\newline
\noindent
  This note is an attempt to redefine event counting problem (cf. \cite{Freq1}, \cite{Freq2}, \cite{Freq3}). In many cases, the most important factor is recent event "popularity rank"  (cf. e.g. \cite{Freq3}) and not its long-run frequency.  Hence, instead of   $n$ item-event counters  consider  a time-dependent discrete probability distribution $ P = (p_1, p_2 \cdots, p_n)  $ as an estimate for relative frequencies (ranks) of the  items involved.   An occurrence of an event with index $i$ can be represented by a delta function distribution $\delta_i $ on the set $ \{1,\cdots,n\} $ triggering an update of estimated probability distribution $P$  by an application of a convex mixture rule $ P \rightarrow  \alpha P + ( 1 - \alpha) \delta_i  $. In other words, arrival of an event $i$ reduces ranks of all other events while tilting estimated  event rank-distribution towards event-item $i$ in a simplest way possible.    
Thus we arrive at the following heavy hitters approximation algorithm 
\begin{algorithm}
		 
		 Fix a number $ \alpha < 1 $ that is close to $1 $. 
	 If  an item $ j \in \{ 1,2, \cdots, n \}  $,  was clicked (event number $j$ did occur) set
 $ p_i \rightarrow \alpha  p_i, \;  i = 1, \cdots n, \; i \neq j $ and set 
	$ p_j  \rightarrow  \alpha p_j + 1-\alpha $ 
\end{algorithm}
One practical problem with the above is that all frequencies (probabilities) are updated simultaneously. There are, however, some advantages:
\begin{itemize}
 \item[(1)] decreasing $\alpha $ gives higher priority  to recent events and vice-versa, increasing  $\alpha$  will bias the ranking towards "idling" event items  
 \item[(2)] and therefore, sensitivity of this ranking scheme to new events can be easily controlled (even at runtime) by adjusting just one parameter 
\end{itemize}

\begin{remark}
	 Suppose that it is desirable that an item should loose half of its rank if it was idle while a list it belongs to  was updated $T$ times.
	 It is quite obvious that  this can be achieved by setting parameter $ \alpha $  to   $  \exp( -\log(2)/T) $. For example, if $T = 10 $ then $ \alpha \approx .93 $ (cf. \cite{rank})
\end{remark}
\noindent Close relationship between Algorithm 1 and Bernoulli convolutions (cf. \cite{ber1}) is a subject of the rest of this paper.

\section{Bernoulli convolutions}
 Suppose that incoming event frequencies follow a fixed discrete distribution  $ Q = (q_1, q_2, \cdots q_n) , \; \sum_{i=1}^n  q_i  = 1   $ and let $ Y_t = (y_{1,t} \;, \cdots , y_{n,t}) $ be a probability distribution vector ($ \sum_{i=1}^n  y_{i,t}  = 1 $ for all $t$) of our (relative) frequency estimates at times $t=0,1, \cdots $ .
 Essentially, Algorithm 1 computes a path of a random walk on a standard $(n-1)$-dimensional simplex $ \sigma^{n-1} \in \mathbb{R}^n $ defined by iterative  rule   
 \begin{equation}
  Y_{t+1} = \alpha Y_t + ( 1 - \alpha ) \delta_i  \; \text{with probability} \; q_i , \; i = 1,2,\cdots, n
 \end{equation}   
where $\delta_i $ is an $i$-th vertex of the simplex  $ \sigma^{n-1}$  or, in other words,  the $i$-th unit vector in standard Eucledean coordinates  in $ \mathbb{R}^n $.
The update rule for the $i$-th coordinate on iteration $t+1$ is   
\begin{eqnarray}
 y_{i, t+1} =
\begin{array}{ll} 
\alpha y_{i,t}  &  \text{ with probability } 1- q_i   \\     
 \alpha y_{i,t} + 1- \alpha  & \text{ with probability } q_i 
\end{array}  \label{1d}
\end{eqnarray}
Let's  fix a coordinate for a while,  omitting the index $i$.
Let $ \xi_m, \; m= 1, \cdots, t $ be random biased Bernoulli variables such that $ 	\mathbb{P}( \xi_m = 0 ) = 1-q $ and $ 	\mathbb{P}( \xi_m = 1 ) = q $.
It is well known (see. e.g. \cite{ber1}) that on step  $t$ the one-dimensional random walk (\ref{1d})   corresponds to a random variable 
\begin{equation}
 y_t = \alpha^t y_0 +   ( 1 - \alpha ) \sum_{m=0}^{t} \xi_m \alpha^m
\end{equation}   
which up to a mostly irrelevant free term is a convolution of $t$ biased Bernoulli variables.  The infinite biased Bernoulli convolution (cf. e.g. \cite{ber1}) is obtained from (3) by setting $t=\infty$ or similarly, by driving the random process (\ref{1d}) infinite number of steps.  
\begin{remark}
 It is well known (see e.g. \cite{shm} for precise statement)  that  Bernoulli convolution $   ( 1 - \alpha ) \sum_{0}^{\infty} \xi_m \alpha^m $      is absolutely continuous (with respect to the Lebesgue measure on the line)  for almost all  sufficiently large values of parameter $\alpha$.  For these values of $\alpha$ the weak limit $y$ of the sequence of random variables $y_t$ does exit and only this case will be considered in this paper. 
\end{remark}
\begin{lemma}
\begin{equation}
	\mathbb{E}(y_t) = \alpha^t y_0 + ( 1 - \alpha^t)q \\
\end{equation}
\end{lemma}
%
Indeed, by definition (2)
\begin{equation}
	\mathbb{E}(y_{t}) = \alpha 	\mathbb{E}(y_{t-1}) + ( 1 - \alpha)q 
\end{equation}   
and hence by  induction
\begin{equation}
	\mathbb{E}(y_{t}) = \alpha^t y_0 +  (1-\alpha) (1 + \alpha + \cdots + \alpha^{t-1}) q  \nonumber
\end{equation}
which is the same as (4).
%
%
\begin{lemma}
	\begin{equation}
	\mathbb{VAR}(y_t) = (1 - \alpha^{2t}) \; \frac{1-\alpha}{1+\alpha} \; ( q - q^2)
	\end{equation}
	\end{lemma}
Proof. It follows from the definition (2) that 
\begin{equation}
	\mathbb{E}(y_{t}^2) = \alpha^{2} 	\mathbb{E}(y_{t-1}^2) +  (1-\alpha)^2 q +
	2 \alpha (1 - \alpha )  \mathbb{E}(y_{t-1}) q   \nonumber
\end{equation}
and therefore by (5) 
\begin{eqnarray}
	\mathbb{VAR}(y_t)  = 	\mathbb{E}(y_{t}^2) - 	\mathbb{E}(y_{t})^2 =  \alpha^2 	\mathbb{VAR}(y_{t-1}) + ( 1 - \alpha)^2 \;  ( q - q^2) 
\end{eqnarray}

\noindent From here,  by  the same inductive argument as in Lemma 1, we get
\begin{equation}
	\mathbb{VAR}(y_t)  = \frac{(1 - \alpha^{2t})} {1-\alpha^2} \; ( 1 - \alpha)^2 \; ( q - q^2 )= (1 - \alpha^{2t}) \; \frac{1-\alpha}{1+\alpha} \; ( q - q^2)   \nonumber
\end{equation}
\noindent 
As an obvious consequence of lemmas 1 and 2 (cf. Remark 2) we have
\begin{corollary}
	The infinite  Bernoulli convolution defined by (2) has  expectation $q$ and variance 
	$ \frac{1-\alpha}{ 1 + \alpha }  ( q - q^2 )  $

\end{corollary}
\begin{remark}
	Under assumption  that the sought for limits exist (Remark 2),  Corollary 1 can be established by passing to the limit in recurrent relations (5), (7) and then solving for expectation and variance respectfully.  
\end{remark}
Here is an example,  demonstrating that passing to a limit as suggested in Corollary 1 is not always possible. 
\begin{example}
Assuming that starting point of the random walk (2) is non-zero, we have 
\begin{eqnarray}
\mathbb{E}\left( \frac{1}{y_{t+1}}  \right) \; = \; ( 1 - q ) \mathbb{E}\left( \frac{1}{ \alpha y_{t}}  \right) \; +
\; q \mathbb{E}\left( \frac{1}{ \alpha y_{t} + 1 - \alpha}  \right)    = \nonumber \\
 = \; \frac{ 1 - q } {\alpha} \; \mathbb{E}\left( \frac{1}{ \alpha y_{t}}  \right) \; + \;
q \mathbb{E}\left( \frac{1}{  1 - \alpha ( 1 - y_{t} ) }  \right)  \nonumber
\end{eqnarray} 
Passing here to the limit as  $ t \rightarrow \infty $ yields 
\begin{equation}
\frac{\alpha - 1 + q }{\alpha} \; \mathbb{E}\left( \frac{1}{y}  \right) \; = \; q \mathbb{E} \left( \frac{1}{  1 - \alpha ( 1 - y ) }  \right)  
\end{equation} 
which is obviously wrong if $ \alpha \leq 1 - q $ and therefore, the condition $ \alpha > 1 - q $ is necessary for   
the existence of continuous limit $ \lim_{t \rightarrow \infty} 1/y_t $. 
 If $ q > 1 - q $ the condition $ \alpha  > 1 - q  $ follows from the well known necessary condition $ \alpha > q^q ( 1-q)^{1-q} $ for non-singularity of Bernoulli convolution $ \lim_{t\rightarrow\infty} y_t $ (cf. e.g. \cite{shm}). For (8) to be true, however, we need non-singularity of the inverse of Bernoulli convolution.
 Essentially a question one can ask is this. For what values of $\alpha$ (if any) $ \lim_{t \rightarrow \infty} 1/y_t $ satisfying (8) exists. 
   
\end{example}

\section{Random walk on a simplex}
 We will compute variances of random vectors  generated by (1) and some other similar random walks.
As before, it is assumed that continuous limit $ Y= \lim_{t\rightarrow \infty} Y_t$
 does exist. It follows from (4-5) and Corollary 1 that     
\begin{eqnarray}
\mathbb{E}( Y_{t+1} )  = \alpha \mathbb{E}( Y_t )  + ( 1- \alpha) Q    \\
\mathbb{E}( Y_t )  = \alpha^t  Y_0 + ( 1- \alpha^t) Q   \nonumber \\
\mathbb{E}( Y)  =  Q   \nonumber
\end{eqnarray}
\noindent In what follows, all vectors are assumed to be column vectors so that for vectors $A,B $ their outer product is $AB^T$ where $ B^T$ is a row vector transposition of $B$. A diagonal matrix with elements of a vector  $ A $  on its  main diagonal will be denoted by    $ diag(A) $. 
\newline
\newline
\noindent  
  Using the rule (1) we get 
\begin{eqnarray}
 \mathbb{E}( Y_{t+1}Y_{t+1}^T ) = \alpha^2 \mathbb{ E}( Y_{t}Y_{t}^T)  + 
( 1 - \alpha)^2 \sum_{i=1}^n q_i  \delta_i \delta_i^{T}   + \nonumber \\
+ \; \alpha(1-\alpha)\sum_{i=1}^n q_i \; (\; \mathbb{E}( Y_t  \delta_i^{T} )  +  \mathbb{E} ( \delta_i Y_t^{T} ) \; ) =  \nonumber \\ = 
 \alpha^2 \mathbb{ E}( Y_{t}Y_{t}^T)  + ( 1 - \alpha)^2 diag(Q) +  \alpha ( 1 - \alpha ) (\; \mathbb{E}( Y_t )Q^T +Q \mathbb{E}( Y_t^{T} )\;)
  \label{long}
\end{eqnarray}
In the same way, using (9) we compute 
\begin{eqnarray}
\mathbb{E}( Y_{t+1} ) \mathbb{E}( Y_{t+1}^T )  = \alpha^2 \mathbb{E}( Y_{t} ) \mathbb{E}( Y_{t}^T )  + ( 1 - \alpha)^2 QQ^T  
+  \alpha ( 1 - \alpha ) (\; \mathbb{E}( Y_t )Q^T +Q \mathbb{E}( Y_t^{T} )\;)   \nonumber
\end{eqnarray}
and subtracting this from (10) we obtain a recurrent relationship  
\begin{equation}
	\mathbb{VAR}( Y_{t+1} )  = \alpha^2 	\mathbb{VAR}( Y_t ) + ( 1 - \alpha)^2   ( diag( Q ) - QQ^T )   \nonumber
\end{equation} 
which is perfectly similar to (7).  Hence, in accordance with Lemma 2 we have 
\begin{theorem}
	The covariance matrix of the finite $n$-dimensional Bernoulli convolution defined by (1) is 
\begin{equation}
		\mathbb{VAR}(Y_t) = ( 1 - \alpha^{2t}) \frac{1-\alpha}{1 + \alpha } ( diag( Q ) - QQ^T)  \nonumber
\end{equation}     
The covariance matrix of the corresponding infinite $n$-dimensional Bernoulli convolution is 
\begin{equation}
	 	\mathbb{VAR}( Y ) = \frac{1-\alpha}{1 + \alpha } \; ( diag( Q ) - QQ^T) \nonumber
\end{equation}  
\end{theorem}
\noindent Let $1_n$ be  $n$-vector with all its coordinates being equal to one. It's easy to check that  
$	\mathbb{VAR}(Y_t)(1_n) =   	\mathbb{VAR}(Y)(1_n)  =  0 $. This is not surprising since coordinates of $ Y_t$  sum-up to one.
The matrix $diag(Q) - QQ^T$ is a symmetric rank-one perturbation of a diagonal  matrix and spectral structure of such matrices is well studied. We just mention
\begin{corollary}
	If  bias probabilities $q_i$ are pairwise distinct then all the non-zero  eigenvalues of the covariance matrix of  $n$-dimensional Bernoulli convolution (1) are   distinct roots of the equation 
\begin{equation}
\sum_{i=1}^n \frac{q_i} {q_i - \lambda} \;  = \;0 \nonumber
\end{equation} 
\end{corollary}
On the other hand, we have
\begin{example}
The only eigenvalues of the covariance matrix of unbiased ($ q_i = 1/n,\; i = 1, \cdots, n$)  $n$-dimensional Bernoulli convolution are $ 0 $ and $ 1/n$
\end{example}
As a slight generalization of (1), fix   $ m > 1 $ points (vectors)
$ v_1, \cdots v_m $ in $ \mathbb{R}^n $ and discrete probability distribution $Q = (q_1, q_2, \cdots q_m ) $. Define a random walk by a rule 
\begin{equation}
 Y'_{t+1} = \alpha Y'_t + ( 1 - \alpha ) v_i  \; \text{with probability} \; q_i , \; i = 1,2,\cdots, m   
\end{equation}  
Let $ V $ be an $ n \times m $ matrix that has coordinates of  $ v_1, \cdots, v_m $ as its columns. For random vectors defined by (11), the equation (9)  turns into
\begin{equation}
\mathbb{E}( Y'_{t+1} )  = \alpha \mathbb{E}( Y'_t )  + ( 1- \alpha) VQ  \nonumber
\end{equation} 
Let $ Y' = \lim_{t\rightarrow \infty} Y'_t $. From the proof of Theorem 1 we have 
\begin{corollary}
\begin{eqnarray}
  \mathbb{E} ( Y'_t ) = \alpha^t Y'_0 + ( 1 - \alpha^t)VQ , \; \;\;   \mathbb{E} ( Y') = VQ           \nonumber \\
 	\mathbb{VAR}(Y'_t) = ( 1 - \alpha^{2t}) \frac{1-\alpha}{1 + \alpha } V( diag(Q ) - QQ^T)V^T  \nonumber \\ 
 		\mathbb{VAR}(Y') =  \frac{1-\alpha}{1 + \alpha } V( diag(Q ) - QQ^T)V^T  \nonumber
 \end{eqnarray}	
\end{corollary}
and in one-dimensional case 
\begin{corollary}
	\begin{eqnarray}
  \mathbb{E} ( Y'_t ) = \alpha^t Y'_0 + ( 1 - \alpha^t)\sum_{i}^m v_i q_i   , \; \;\;   \mathbb{E} ( Y') = \sum_{i}^m v_i q_i            \nonumber \\
\mathbb{VAR}(Y'_t) = ( 1 - \alpha^{2t}) \frac{1-\alpha}{1 + \alpha } \sum_{i=1}^n v_i^2 ( q_i - q_i^2)  \nonumber \\ 
\mathbb{VAR}(Y') =  \frac{1-\alpha}{1 + \alpha } \sum_{i=1}^n v_i^2 ( q_i - q_i^2)    \nonumber
\end{eqnarray}
\end{corollary} 
Note that setting here  $ m = 2,  v_1 = 0, v_2 = 1 $ we not-surprisingly recover equations  (4) and (6).
\newline
\newline
\noindent 
Moreover, consider a case when all points $ v_1, \cdots , v_m $ belong to a complex plain. Then $ Y'_t, \; t = 1, 2, \cdots  $ is  a sequence of complex random variables and  again from the proof of Theorem 1 we have 
\begin{theorem}
	Let $ v_1, \cdots, v_m \in \mathbb{C}^1, \; m > 1 $. Then for the sequence of complex random variables $ Y'_t, \;  t = 1, \cdots  $ defined by (11) we have  
	\begin{eqnarray}
\mathbb{E} ( Y'_t ) = \alpha^t Y'_0 + ( 1 - \alpha^t)\sum_{i}^m v_i q_i   , \; \;\;   \mathbb{E} ( Y') = \sum_{i}^m v_i q_i            \\
\mathbb{VAR}(Y'_t) = ( 1 - \alpha^{2t}) \frac{1-\alpha}{1 + \alpha } \left(\sum_{i=1}^n |v_i|^2 ( q_i - q_i^2) - \sum_{i < j} ( v_i \bar{v}_j + v_j \bar{v}_i) q_i q_j   \right) \nonumber \\ 
\mathbb{VAR}(Y'_t) =\frac{1-\alpha}{1 + \alpha } \left( \sum_{i=1}^n |v_i|^2 ( q_i - q_i^2) - \sum_{i < j} ( v_i \bar{v}_j + v_j \bar{v}_i) q_i q_j   \right) \nonumber 
	\end{eqnarray} 
\end{theorem} 
The proof is similar to the proof of Theorem 1. By definition
\begin{equation} 
\mathbb{VAR}(Y'_t)  =  \mathbb{ E} (Y'_t \bar{Y'_t})  ) - \mathbb{E} (Y'_t) \mathbb{ E}( \bar{Y'_t})     \nonumber
 \end{equation} 
and as in the proof of Theorem 1
\begin{equation}
\mathbb{E}( Y'_{t+1} \bar{Y'}_{t+1})   = \alpha^2  \mathbb{E}(Y'_{t}\bar{Y'}_{t} )  + 2 \alpha( 1- \alpha) 
\mathbb{E}(Y'_t)  \mathbb{E}(\bar{Y'}) +   ( 1 - \alpha )^2   \sum_{i=1}^n q_i |v|_i^{2} 
\end{equation}
On the other hand
\begin{equation} 
\mathbb{E}( Y'_{t+1} )  \mathbb{E}(\bar{Y'}_{t+1})   = \alpha^2  \mathbb{E}(Y'_{t} )  \mathbb{E} ( \bar{Y'}_{t} )  + 2 \alpha ( 1- \alpha) 
\mathbb{E}(Y'_t)  \mathbb{E}(\bar{Y'}) +   ( 1 - \alpha )^2 \mathbb{E}(\bar{Y'})  \mathbb{E}(Y')  \nonumber
\end{equation}
and it follows from (12) that 
\begin{equation}
 \mathbb{E}(\bar{Y'})  \mathbb{E}(Y') =  \sum_{i=1}^n |v_i|^2  q_i^2 +  \sum_{i < j} ( v_i \bar{v}_j + v_j \bar{v}_i) q_i q_j    \nonumber
\end{equation}
Substituting this into previous equation and subtracting from (13) we obtain a recurrent relation  
\begin{equation}
\mathbb{VAR} (  Y'_{t+1} ) = \alpha^2 \mathbb{VAR} (  Y'_{t} )  + ( 1 - \alpha)^2 \left( \sum_{i=1}^n |v_i|^2 ( q_i - q_i^2) - \sum_{i < j} ( v_i \bar{v}_j + v_j \bar{v}_i) q_i q_j   \right)  \nonumber
\end{equation}
The rest of the proof is the same as in  Lemma 2.

\begin{corollary}
	If all points $ v_i \in \mathbb{C}, \; i = 1,2, \cdots, m, \; m > 1 $  belong to a unit circle then
\begin{equation}
\mathbb{VAR}(Y'_t) = 4 ( 1 - \alpha^{2t})  \frac{1-\alpha}{1 + \alpha }  \sum_{i < j}  \sin^2 \left( \frac{ \phi_{i,j}} {2} \right) q_i q_j       \nonumber \\ 
\end{equation}	 
where   $ \phi_{i,j},  \; i < j $  are pairwise angles between unit vectors  $v_i, v_j $. 
\end{corollary}
Indeed, since in this case  $ |v_i | = 1, \; i = 1, \cdots m $ , we have 
\begin{equation}
 \sum_{i=1}^m |v_i|^2  q_i = 1,  \;  v_i \bar{v}_j + v_j \bar{v}_i = 2\cos( \phi_{i,j})    \nonumber
\end{equation} 
and on  the other hand
\begin{eqnarray}
\sum_{i=1}^m q_i^2 \; + \;  2 \sum_{i<j} \cos(\phi_{i,j}) q_i q_j  \; =  \ 1  \; +  \; 2 \sum_{i , j} ( \cos(\phi_{i,j} ) - 1  ) q_i q_j  \nonumber \\ 
= 1 - 4 \sum_{i , j}  \sin^2( \phi_{i,j }/2 ) q_i q_j 
\nonumber
\end{eqnarray}
\begin{example}
If $m=n=3$ then  two-dimensional random walk (1) can be viewed as a random walk on an equilateral triangle  $\sigma^2$  whose vertices are three distinct cubic roots of unity  $ v_1 = 1, 
v_2 = e^{2\pi i /3 } , v_3 = e^{4\pi i / 3 } $. All three angles between $ v_i $ and $v_j, \; i,j = 1,2,3, \;
 i \neq j $ are equal to $ 2 \pi /3 $ and  
by Corollary 5 the (complex) variance of the corresponding complex random variable at iteration  $t$ is 
\begin{equation}
\mathbb{VAR}(X_t)  =   3 \; ( 1 - \alpha^{2t})  
\;  \frac{1-\alpha}{1 + \alpha} \; ( q_1q_2 + q_1 q_3 + q_2 q_3  )    \nonumber
\end{equation}
\end{example}
	
\section{Properties of  approximation}	
Results of the  section 2 can be used to evaluate heavy hitters approximation produced by Algorithm 1. 
\newline
\newline
\noindent
To evaluate the algorithm ability to "overweight"  recent event frequencies, let's assume that the number of iterations $  t $  corresponds  to a "relevancy" time window. For example, if last week heavy hitters are of highest importance, let  $t$ be  a "weekfull of clicks".  Measuring time by click-counter, suppose that estimated click-distribution at the start of the time period was  $X$ and that for time $ t_1 $ the incoming click distribution $ P_1 $ did not change. Suppose also that at time $t_1$    the incoming distribution switched to $ P_2  $ and did not change for the remaining time  $ t_2 = t  -   t_1  $.
Then by Lemma 1,  an expected convex mixture approximation at the end of the time period will be 
\begin{equation}
\alpha^{t} X + \alpha^{t_2}(1 - \alpha^{t_1}) P_1 + ( 1 - \alpha^{t_2}) P_2  \nonumber
\end{equation} 
To see how our approximation is affected by recent events  let's estimate  the ratio of coefficients at $P_2$ and $ P_1 $ in the expression above.
Since $ \beta = 1 - \alpha $ is supposed to be small, we have  
\begin{equation}
\frac  {1 - \alpha^{t_2} } {   \alpha^{t_2}(   1 - \alpha^{t_1}   )  }     = 
\alpha^{-t_2}   \frac   { 1 - (1 - \beta)^{t_2}   }  { 1 - ( 1 - \beta )^{t_1}  } \approx \alpha^{-t_2}  \frac{t_2}{t_1}
\end{equation} 


\noindent In case of plain event counting this ratio should be $ \approx t_2/t_1 $. On the other hand,   from (14) we have 
\begin{corollary}
	Algorithm 1 introduces approximately times $\alpha^{-1} $ per iteration "velocity boost" for recent heavy hitters. 
\end{corollary}

\noindent As we saw above, Algorithm 1 will approximate the mean of a  fixed  incoming click distribution  in the long run.  Lemmas 1, 2 and a straightforward application of Chebyshev inequality (cf. e.g. \cite{Fer} for a  vector version) give a reasonable estimate for a quality of this approximation.
\begin{corollary}
The following estimates hold for random variables  $y_i= \lim_{t\rightarrow\infty} y_{i,t} $ and  for random vector $	 Y= \lim_{t\rightarrow\infty} Y_t  $
\begin{equation}
 \mathbb{P} \left( \; | y_{i} - q_i | \; \geq \;  \epsilon  \; \right)   \;  \leq \;  \frac{1 - \alpha }{ \epsilon^2(1 + \alpha)} ( q_i - q_i^{2} ), \; i = 1, \cdots n  
\end{equation}
In particular,
\begin{equation} 
\mathbb{P} \left( \; | y_{i} - q_i | \; \geq \;  \sqrt{1-\alpha}  \; \right)   \;  \leq \;  \frac{q_i - q_i^{2}}{1 + \alpha},    \; i = 1, \cdots n 
\end{equation}
and 
\begin{equation}
	\mathbb{P} \left( \; \parallel Y - Q \parallel \; \geq \;  \epsilon 
		  \; \right)   \;  \leq \;    \frac{1-\alpha}{\epsilon^2(1+\alpha)}  (1 -   \sum_{i} q_i^2 )   \nonumber
\end{equation}
\end{corollary}

\begin{remark}
	It follows from (16) that for any $q_i$ and large enough $\alpha$, about $(7/8)$-th of the  limit distribution belongs to the narrow interval $ [ -\sqrt{1-\alpha}, \sqrt{1-\alpha}]   $
\end{remark}

\begin{example}
	For  $ \alpha = 0.99, \;  q = 1/2 $ and for sufficiently large $t$ the  value of $ y_t $ will belong to the interval $[ 0.4, \; 0.5 ]$ with about  $87\%$ probability
\end{example}

\noindent 
 It is obvious, that the estimator  (15)  works better for large values of $q$, i.e. for above-mentioned heavy hitters.   More precisely, setting $ \epsilon \leftarrow \epsilon q_i $ in (15) we get
\begin{corollary}
	An estimate
		\begin{equation}  
	\mathbb{P} \left(  \; | y_i - q_i| \; \geq \epsilon q_i     \; \right)   \;  \leq \; \epsilon   \nonumber 
	\end{equation}
	holds  for 
\begin{equation}
q_i  \; \geq \;  \frac{ 1 }  {    1 + \frac{ 1+ \alpha }  {1-\alpha }   \epsilon^3    }    \nonumber
\end{equation}  

\end{corollary}

\begin{example}
  For $ \epsilon = 1/10 $ and $ \alpha = 1 - \epsilon^3 = .999 $ this boils down to  
	\begin{equation}  
	\mathbb{P} \left(  \; | y_i - q_i | \; \geq \;   q_i/10     \; \right)   \;  \leq \; 1/10   
\;\;	\text{ if } \;\;	q_i \; \geq \; \;  \frac{ 1 }  {    2.999     }     \nonumber	
	 \nonumber
	\end{equation}
In other words, for large enough number of iterations,  click probabilities that are slightly above $1/3$  can be  approximated  up-to  $10\% $ relative error with $90\%$ confidence.
\end{example}
\noindent For a finite Bernoulli convolutoin obtained after $t$ iterations of Algorithm 1 we  get from (4) and (6)
\begin{corollary}
	If $y_{i,0} = q_i, \; i = 1, \cdots n $ then for any $ t = 1,2, \cdots $ 
	\begin{equation} 
	\mathbb{P} (\; | y_{i,t} - q_i | \; \geq  \epsilon \;)   \;  \leq \;  \frac{ (1 - \alpha^{2t} ) (1-\alpha ) }{  \epsilon^2(1 + \alpha) }  \; ( q_i - q_i{^2} ),     \; i = 1, \cdots n \nonumber 
	\end{equation}
	In particular 
		\begin{equation} 
	\mathbb{P} (\; | y_{i,t} - q_i | \; \geq  \sqrt{1-\alpha} \;)   \;  \leq \;  \frac{ (1 - \alpha^{2t} )  }{ 1 + \alpha }  \; ( q_i - q_i{^2} ),     \; i = 1, \cdots n \nonumber 
	\end{equation}
	and if $ Y_0 = Q$ then
	\begin{equation}
	\mathbb{P} \left( \; \parallel Y_t - Q \parallel \; \geq \;  \epsilon 
	\; \right)   \;  \leq \;    \frac{  ( 1 - \alpha^{2t} ) ( 1-\alpha)}{\epsilon^2(1+\alpha)}  (1 -   \sum_{i} q_i^2 )   \nonumber
	\end{equation}
\end{corollary}

\section{Recurrent formula for moments  of biased Bernoulli convolutions}
    
  Moments of unbiased Bernoulli convolutions were studied in  \cite{moments1},\cite{moments2},\cite{moments3}. 
Some basic properties of  moments of biased infinite Bernoulli convolutions are briefly discussed in this section.. 
\newline
\newline
\noindent  It makes sense to consider central moments, 
  $ \mathbb{E}(y  -  q)^n $  (cf. Corollary 1). Hence, we replace the sequence 
  $ y_t $ with the sequence $ y_t - q $ which from now on will be denoted by the same letter. The transformation rule (2) thus changes to 
\begin{eqnarray}
y_{m+1} =
\begin{array}{ll}
\alpha y_{m}  - ( 1- \alpha ) q &  \text{ with probability } 1- q   \\     
	\alpha y_{m} + (1- \alpha )(1 - q )  & \text{ with probability } q 
	\end{array} 
\end{eqnarray}  
For expectations of the random variable sequence $ y_{m}^n, \; m=1,2,\cdots \; $ that tarnslates into  
\begin{equation}
\mathbb{E}(y_{m+1}^{n} )  = ( 1-q) \mathbb{E}( (\alpha y_{m}  - ( 1- \alpha ) q  )^n) + q\mathbb{ E}(	( \alpha y_{m} + (1- \alpha )(1 - q )  )^n ) \nonumber
\end{equation}
Opening brackets and passing to the limit (that is assumed to exist) results in identity  
\begin{equation}
\mathbb{E}(y^{n} )  = \alpha^n \mathbb{E}( y^n)   +   \sum_{k=1}^n \binom{n}{k} \alpha^{n-k}( 1- \alpha )^{k}( ( -q)^k (1-q) + q(1-q)^k) )  \mathbb{E}(y^{n-k})   \nonumber
\end{equation}
Finally,  after relabeling $ M_k = \mathbb{E}( y^k ) $ we obtain for $n>= 2 $ a recurrent relation (cf. \cite{moments2}) 
\begin{eqnarray}
M_n = \frac{1}{1 - \alpha^n}    \sum_{k=1}^n \binom{n}{k} \alpha^{n-k}( 1- \alpha )^{k}( (1-q)( -q)^k + q(1-q)^k) ) M_{n-k}  \nonumber \\
\equiv \; \frac{q-q^2}{1 - \alpha^n}    \sum_{k=2}^n \binom{n}{k} \alpha^{n-k}
( 1- \alpha )^{k}
( ( -1)^k q^{k-1} + (1-q)^{k-1}) ) M_{n-k}  \label{binom}
\end{eqnarray}
Obviously,  $ M_0 = 1$ and $ M_1 = 0 $. It is now a simple matter to write down a few central moments of the infinite Bernoulli convolution (2):

\begin{example}
	\begin{eqnarray}
		M_2 =  \frac{1-\alpha}{ 1 + \alpha }  ( q - q^2 )   \nonumber \\
		M_3 =  \frac{(1-\alpha)^3}{ 1 - \alpha^3 }  ( q - q^2 )(1-2q)   \nonumber \\
	M_4 = \frac{ (1-\alpha)^4 } {1- \alpha^4} ( q - q^2 ) \left[ \frac{6 \alpha^2}{ 1 - \alpha^2} ( q - q^2 )  +  1 -q + q^2  \right]   \nonumber
	\end{eqnarray}
\end{example}
Let  $ \mu_y $ be a measure associated with the infinite Bernoulli convolution  $ y $ that is generated by rule (17)
and let $ (.)^{*} $ denote a reflection  $ x \rightarrow 1 - x $. Denote also by $y^{*}$ an infinite Bernoulli convolution  generated by the  rule (17) with interchanged  probabilities $ q \rightarrow 1 - q $.
It is probably worth mentioning
\begin{corollary}.
\begin{itemize} 

	\item[(i)] for any interval $[a,b], \; \mu_{y^{*}}([a,b]^{*}) =  \mu_y ( [a,b]) $ 
	\item[(ii)] $ y^{*} = -y $ and therefore
	\item[(iii)] $\mathbb{E}( y^{n} ) = (-1)^n \; \mathbb{E}( y^{*n} ),  \;\; n = 0,1,2 \cdots $
	\item[(iv)] as polynomials of $q$, the central moments $ M_n(q) \equiv \; \mathbb{E}(y^{n})(q)  $ are semi-invariant with respect to the involution $ \tau : q \rightarrow 1- q $, that is
\begin{equation}
M_n^{\tau}(q) = (-1)^n M_n(q)  \nonumber
\end{equation}
\end{itemize} 
\end{corollary}
Indeed, statements (i) and (ii) follow from definition (17).
Statement (iii)  follows from (ii) or (iv) and  the proof of (iv) is a straightforward induction based on (18). 
\newline
\newline
\noindent Moreover, for  central moments $M_n \equiv M_n(q) $ as polynomials of $q$ we have  
\begin{corollary}
 $ M_n(q)$ is a polynomial of $ q (1 - q ) $ if n is even and is a polynomial of $ q(1-q) $ times $ 1 - 2q $ if $ n $ is odd.  
\end{corollary}
This is an easy consequence of Corollary 10. Just note, that  it follows from Corollary 10 (iv) that  $M_n(q)$ is divisible by $ q - \frac{1}{2} $ if $n$ is odd.
\begin{lemma}
	If $ q \leq 1 - q $ then  
\begin{itemize}
	\item[(i)] 	all central moments $ M_n$ are non-negative
	\item[(ii)]  $ M_n \leq 1 -q $ for all $ n = 0,1, \cdots $
	\item[(iii)] $ \lim_{n\rightarrow \infty}M_{2n}^{1/2n} = 1 - q $
\end{itemize}
\end{lemma}	
Proof. The first statement  directly follows from (18). The second statement is 
obvious. Statement (iii) is just a recollection of a well known fact about a sequence of $n$-norms $  \left( \int_{-q}^{1-q} |y(x)|^n d\mu_y(x) \right)^{1/n} $ converging to $ \infty$-norm $ \max \{ |y| \} = 1 -q $.
\newline
\newline
\noindent Although random variable $ y $ is not non-negative, the following still holds
\begin{theorem}
	If $ q \leq 1 - q $ then $ \lim_{n\rightarrow \infty}M_n^{1/n} = 1 - q $
\end{theorem}
Proof. The sequence $ M_n^{1/n}, \; n = 0, 2, \cdots$ for even numbered central moments is non-decreasing  by H\"older's  inequality and converges to $1-q$  by Lemma 3. 
Hence, for any $ \epsilon_1 > 0 $ there is $k= k_0 $ such that 
\begin{equation}
 M_n \geq ( 1 - q - \epsilon_1 )^n 
\end{equation}
 for all even $ n $ such that $ n \geq  k_0 $. In particular  $ M_{n-k} \geq ( 1 - q - \epsilon_1 )^{n-k} $
for all odd $ k $ and $ n $  such that  $ k \leq  n - k_0 $.  
Using this fact, we will show that an estimate similar to  (19) holds for any large odd number  $n$. Indeed, 
it follows  from (18), Lemma 3 (i) and (19) that for any odd $ n > k_0 + 2 $  
\begin{equation}
M_n \; \geq \;  \frac{q( 1 - q - \epsilon_1 )^n}{1- \alpha^n} \sum_{ 2 \leq k \leq n - k_0,\; k \; \text{odd} } 
\binom{n}{k}
 \alpha^{n-k}
 ( 1 - \alpha)^k          
\end{equation}
It is easy to see, however,  that the sum in (20) is equal to  
\begin{equation}
 \frac{1}{2} \; - \;  \alpha^n \; - \sum_{ n- k_0 < k \leq n,\; k \; \text{odd} } 
\binom{n}{k}\alpha^{n-k}( 1 - \alpha)^k   \nonumber
\end{equation}
and therefore for any $ \epsilon_2 > 0 $ we can find large enough $ n_0 $ such that for any odd  
$ n > n_0 $
\begin{equation}
\sum_{ 2 \leq k \leq n - k_0,\; k \; \text{odd} } 
\binom{n}{k}\alpha^{n-k}( 1 - \alpha)^k  \; \geq \;  \frac{1}{2} - \epsilon_2  \nonumber
\end{equation} 
After substituting this into (20) we find that
\begin{equation}
M_n^{1/n}  \; \geq \; ( 1 - q - \epsilon_1 ) \left(  \frac{q}{1- \alpha^n} \right)^{1/n} 
\left( \frac{1}{2}  - \epsilon_2 \right)^{1/n}    \nonumber
\end{equation}
which is a desired estimate of $ M_n$   for large enough odd $ n $.

\section{Concluding remarks}
As was shown above, relative heavy hitters can be approximated by iterative application of the  convex mixture rule (1). Suggested  algorithm  essentially computes  a Bernoulli convolution if and while an incoming click distribution remains fixed. In practice, the stochastic process of incoming events is much more complicated (cf. e.g. Corollary 3). A problem  of  obtaining similar  convex mixture approximation estimates in a general setting of varying incoming click distributions seems to be both hard and interesting.

\end{document}